


\input amstex

\def\bcc#1{\Bbb C^{#1}}
\def\brr#1{\Bbb R^{#1}}

\def\cf{\Cal F}

\def\ee#1{e_{#1}}
\def\dee#1#2{e_{#1#2}}

\def\ff#1{\Bbb F\Bbb F^{#1}}
\def\fff#1#2{\Bbb F\Bbb F^{#1}_{#2}}

\def\ii{| II |}

\def\La#1{\Lambda^{#1}}
\def\na{n+a}

\def\ooo#1#2{\omega^{#1}_{#2}}
\def\oo#1{\omega^{#1}_0}

\def\ot{\!\otimes\!}

\def\pp#1{\Bbb P^{#1}}
\def\ppp{\Bbb P}
\def\qq#1#2#3{q^{#1}_{{#2} {#3}}}

\def\ra{\rightarrow}

\def\rr#1#2#3#4{r^{#1}_{{#2} {#3}{#4}}}

\def\sx{\sigma (X)}
\def\tx{\tau (X)}

\def\tann{\text{Ann}\,}
\def\tdim{\text{dim}\,}

\def\tsingloc{\text{singloc}\,}
\def\tbaseloc{\text{Base}\,}

\def\up#1{{}^{({#1})}}

\def\ww{\wedge}

\documentstyle{amsppt}
\magnification = 1200
\hsize =15truecm
\hcorrection{.5truein}
\baselineskip =18truept
\vsize =22truecm
\NoBlackBoxes
\topmatter
\title On the infinitesimal rigidity of homogeneous varieties
\endtitle
\rightheadtext{rigidity of homogeneous varieties}
\author
  J.M. Landsberg
\endauthor

\abstract{Let $X\subset\pp N$ be a variety (respectively a
patch of an analytic submanifold) and let $x\in X$ be a general point.
We show that if the projective second fundamental form of
$X$ at $x$ is isomorphic to the second fundamental form
of a point of a Segre $\pp n\times\pp m$, $n,m\geq 2$,   a
Grassmaniann $G(2,n+2)$, $n\geq 4$,
or the Cayley plane $\Bbb O\pp 2$, then $X$ is the corresponding
homogeneous variety (resp. a patch of the corresponding
homogeneous variety). If the projective second fundamental form of
$X$ at $x$ is isomorphic to the second fundamental form
of a point of a Veronese $v_2(\pp n)$ and the 
Fubini cubic form
of $X$ at $x$ is zero, then $X=v_2(\pp n)$ (resp. a patch of
$v_2(\pp n)$).
All these results are valid in
the real or complex analytic categories and 
locally in the $C^{\infty}$
category if one assumes the hypotheses hold in a neighborhood
of any point $x$. As a byproduct, we show that  the systems of
quadrics
$I_2(\pp {m-1}\sqcup\pp {n-1} )\subset
S^2\bcc{m+n}, I_2(\pp 1\times\pp{n-1})\subset S^2\bcc{2n}$
and $I_2(\Bbb S_5)\subset S^2\bcc{16}$
are stable  
in the sense
that if $A_t\subset S^2T^*$ is an analytic family  
such that for $t\neq 0$, $A_t\simeq A$, then
$A_0\simeq A$.

 We also make some observations related to
the Fulton-Hansen connectedness theorem.}
\endabstract

\date {October 26, 1997}\enddate
\address{ Laboratoire de Math\'ematiques,
  Universit\'e Paul Sabatier, UFR-MIG
  31062 Toulouse Cedex 4,
  FRANCE}\endaddress
\email {jml\@picard.ups-tlse.fr }
\endemail
\thanks {Supported by  NSF grant DMS-9303704.}
\endthanks
\endkeywords
\subjclass{ primary 14, secondary 53}\endsubjclass
 
\endtopmatter

\document

The intrinsic rigidity of homogenous spaces has
been studied extensively (see  [HM] and the
references therein). In this paper
we study the extrinsic and infinitesimal rigidity of homogeneous spaces.

Let  $G/P\subset\pp N$ be an
$n$-dimensional homogenous space embedded homogeneously,
but not necessarily in its canonical embedding. Let
$X^n\subset\pp\na$ be a variety and let $x\in X$ be a general point
(that is a point where all integer valued differential invariants
are locally constant).
This paper addresses the question:
To what extent do the  projective differential invariants of $X$
at $x$ need to resemble those of a point of $G/P\subset\pp N$
in order to
be able to  conclude $X=G/P$ as a projective variety?
Let $|II_{X,x}|\subset\ppp S^2T^*_xX$ denote (the system of quadrics
induced by) the projective second fundamental form of $X$ at $x$.
Our progress on this question is as follows:

\proclaim{Theorem 1}Let $X^{n+m}\subset\pp{nm+n+m+z}$, $n,m\geq 2$, be a
patch of an analytic manifold not contained
in a hyperplane and let
$x\in X$ be a general point. If the second fundamental
form $|II_{X,x}|$ is isomorphic to 
$I_2(\pp {m-1}\sqcup\pp {n-1} )$ (the quadrics vanishing on the disjoint
union of two projective spaces), then $z=0$
and $X$ is an open
subset of the Segre $\pp n\times \pp m\subset\pp{nm+n+m }$.
The same result holds locally in the $C^{\infty}$ category if
  the hypotheses hold in the neighborhood of any point
$x$.
\endproclaim

\proclaim{Theorem 2}$X^{2(m-2)}\subset\pp{\binom m2 -1+z}$, $m\geq 6$, be a
patch of an analytic manifold not contained
in a hyperplane and let
$x\in X$ be a general point. If the second fundamental
form $|II_{X,x}|$ is isomorphic to
$ I_2(\pp 1\times\pp{n-1}) $
(the quadrics vanishing on the Segre variety), then
$z=0$ and $X$ is an open
subset of the Grassmanian $G(2,m)$.
The same result holds locally in the $C^{\infty}$ category if
  the hypotheses hold in the neighborhood of any point
$x$.
\endproclaim

Note that the result is false for $m=4$.

\proclaim{Theorem 3}Let $X^{16}\subset\pp{26+z}$ 
be a patch of an analytic manifold not contained
in a hyperplane 
and let
$x\in X$ be a general point. If the second fundamental
form $|II_{X,x}|$ is isomorphic to
$I_2(\Bbb S_5) $
(the quadrics vanishing on the spinor variety), then
$z=0$ and $X$ is an open
subset of the Cayley plane in its canonical embedding.
The same result holds locally in the $C^{\infty}$ category if
  the hypotheses hold in the neighborhood of any point
$x$.
\endproclaim

 The {\it Fubini cubic form} of $X$ at $x$, $F_{3X,x}$ is a relative
  differential invariant encoding the geometric information in
the third derivatives of the embedding. It was
  first used
by Fubini [F] to study hypersurfaces. See [L1] for a precise definition.

\proclaim{Theorem 4} Let $X^n\subset\pp{\binom{n+2}2-1}$,
$n>1$, be a patch of an
analytic manifold not contained in a hyperplane and let $x\in X$ be
a general point. If $III_{X,x}=0$ and $F_{3X,x}=0$, then $X$ is
an open subset of the Veronese $v_2(\pp n)$.
The same result holds locally in the $C^{\infty}$ category if
  the hypotheses hold in the neighborhood of any point
$x$.
\endproclaim

\proclaim{Theorem 5} The systems of quadrics, 
$A= I_2(\pp {m-1}\sqcup\pp {n-1} )\subset
S^2\bcc{m+n},I_2(\pp 1\times\pp{n-1})\subset S^2\bcc{2n}$
and $I_2(\Bbb S_5)\subset S^2\bcc{16}$ are stable in the sense
that if $A_t\subset S^2T^*$ is an analytic family  
such that for $t\neq 0$, $A_t\simeq A$, then
$A_0\simeq A$.\endproclaim

Note that in the analytic category, 
the results are also global, since a patch determines an
entire variety. One must
be careful in the case of the real Cayley plane to insure that
the normalization of the
second fundamental form used in computations is 
possible over $\brr{}$.

\subheading{Previous results}  Monge showed that  a curve in $\pp 2$ is a
conic  
if and only if a   fifth order invariant is zero (see   [L1, 3.6]). In
higher dimensions, Fubini  showed that to determine 
if
a hypersurface is a quadric,   all third order invariants must be
zero [F]. In another direction, it is known
that the Segre variety cannot be deformed as a
submanifold of  projective space (see, e.g. [HM]
section 3).  Note that while $\pp 1\times\pp 1$ 
is rigid as a submanifold of
projective space in the sense of [HM],
it fails to satisfy the analog of theorem 1.
 In [L2] we showed that to determine if $X$ is one of the four
Severi varieties
(that is, $\Bbb A\pp 2$ in its canonical embedding,
where $\Bbb A$ is the complexification of one of the four
real division algebras, i.e.
$v_2(\pp 2)\subset\pp 5, Seg(\pp 2\times\pp 2)\subset\pp 8,
G(2,6)\subset\pp {14}, E_6/P_1\subset\pp{26}$), 
it is necessary to have
agreement of second fundamental forms and a partial vanishing of the cubic
form.
(The case of $v_2(\pp 2)$ had been proven earlier by Griffiths and Harris
[GH].) Theorems 1-3 strengthen these results for the three Severi
varieties with degenerate tangential varieties
in the sense  that they show it is only necessary 
to have  agreement
of second fundamental forms. This strengthening follows
immediately from proposition 7 below.

In the Euclidean geometry of submanifolds, if the Euclidean second
fundamental form is surjective, then a submanifold is
uniquely determined by second order data (sometimes
even first, e.g. hypersurfaces of large dimension). In projective
geometry, the
order of data   needed
to obtain a complete set of functionally independent differential
invariants
  is not    known
except in some special cases. For curves in $\pp 2$ 
 sixth order information is necessary and sufficient.
For hypersurfaces of dimension greater
than two, Jensen and Musso proved third order information is
necessary and sufficient
[JM].  

\subheading{Intrinsic vs extrinsic geometry}
The intrinsic and extrinsic
geometries of homogeous spaces are closely related.
Given any $G/P$,
the
cone of minimal degree rational curves passing through
each point   is an
essentially intrinsic object. The projectivization
of this cone is the base locus of
$|II|$ in the canonical embedding.
(A line osculating to order two at a point
of a homogeneous space  $G/P$ is contained in $G/P$
since homogeneous varieties are cut out by quadrics.)
 
   The intrinsic rigidity  results in  [HM] resemble
ours and have a similar method of proof. 
Hwang and Mok prove
K\"ahler rigidity of Hermitian symmetric spaces of the compact type under
K\"ahler deformations by studing deformations of the cone of minimal degree
rational curves.  This cone naturally sits in a projective space, 
thus their
study at the infinitesimal level is similar to our extrinsic problem.
However, their   results are different, which can most easily
be seen by the fact that the quadric hypersurface is not rigid
to order two, but is K\"ahler rigid, and even holomorphically
rigid (see [H]).  
 It
would be desirable to rephrase the proofs here in terms of
a geometric property of the base locus of the second fundamental
form as in [HM] (see below).

\subheading{Secant and dual varieties}Extremal
degeneracies of auxilliary varieties 
often force  homogeneity.  
Zak proved that if $X^n\subset\pp\na$ is a smooth
variety not contained in a hyperplane, and $a<\frac n2 +2$
then the secant variety $\sx$ must be equal to $\pp\na$
and if $a=\frac n2+2$ and $\sx\neq\pp\na$, then $X$ must
be a Severi variety (see [Z]). Zak also proved 
  that if $X^n\subset\pp\na=\ppp V$ is a smooth
variety not contained in a hyperplane, then $\tdim X^*\geq\tdim X$,
where $X^*\subset\ppp V^*$ denotes the dual variety of $X$. 
Ein showed that if $\tdim X^* =\tdim X$, and $a\geq\frac n2$,
then $X$ is either a hypersurface, $Seg (\pp 1\times\pp m )$,
the Grassmanian $G(2,5)$ or the ten dimensional spinor
variety $\Bbb S_5$ (all in their canonical embeddings),
see [E].

If $X^n\subset\pp\na$ is
a smooth variety with degenerate secant variety, then
$a\leq\binom{n+1}2$
(see [Z], [L2]). A special case of Zak's theorem on
Scorza varieties states that if $a=\binom{n+1}2$
and $\sigma (X)$ is degenerate, then $X$ must
be a Veronese  $v_2(\pp
n)$.

\subheading{The refined third fundamental form and connectedness}
Let $X^n\subset\pp\na$ be a 
patch of a complex manifold. Let $x\in X$ be a general point
and let $v\in T_xX$ be a generic tangent vector. If the mapping
$II_v: T_xX\ra N_xX$, 
defined by $w\mapsto II(v,w)$, is not surjective,
there is a well defined third order invariant, called {\it the
third fundamental form refined with respect to $v$}, $III^v$
(see [L2] for
details).  Given a system of
quadrics $A\subset S^2T^*$ and $v\in T$, let
$\tann (v)=\{ q\in A\mid [v]\in q_{sing}\}$, and let
$\text{Singloc}(A)=\{ v\in T\mid [v]\in q_{sing}
\,\forall q\in A \}$.
Note that $|II_{X,x}|/\tann (v)$ is a well defined system
of quadrics on $ \text{Singloc}(\tann (v))$.
With these notations, 
$$
III^v\in S^3(\text{Singloc}(\tann (v))^* \ot N_xX/II_v(T).
$$

A
special case of the Fulton-Hansen connectedness theorem [FH]
states
that  if $\tdim\sx\neq 2n+1$ or $\tdim \tx\neq 2n$, 
then $\sx=\tx$ for any projective variety $X$. 

A consequence of the Fulton-Hansen theorem is   that if $X$ 
is a smooth variety with degenerate secant variety, then the
refined third fundamental form is zero at general points.
In fact, the refined third fundamental
form being zero implies $\sigma (X) = \tau (X)$ in the
case
$X$ is smooth, see [L2].
In our original proof of Zak's theorem, we used
the consequence of the connectedness theorem
that $III^v\equiv 0$
to prove the rigidity of varieties that infinitesimally
looked like   Severi varieties.

If $A\subset S^2T^*$ is a system of quadrics, its
{\it prolongation}, $A\up 1$ is defined by
$A\up 1 = S^2T^*\cap ( A\ot T^*)$.

\proclaim{Proposition 6}
Let $X^n\subset\pp\na$ be a variety. Let $x\in X$ be a general 
point and let $v\in T_xX$ be a generic tangent vector. With the
notation of the paragraphs above, consider
$|II_{X,x}|/\tann (v)$ as a system of quadrics on
$\tsingloc (\tann (v))$. Then
$$
|III^v|\subseteq (|II_{X,x}|/\tann (v))\up 1.
$$
\endproclaim

Proposition 7 follows   from the
 formula [L2, 13.1]. The first line of [L2, 13.1]
shows that
$|III^v|\subseteq (|II_{X,x}|/\tann (v))\ot T^*$, and the second
line shows
that it is symmetric.

Zak's theorem on Scorza varieties
 indicates that perhaps theorem 4 is   not
the optimal result.
A positive answer to the following question
 would provide
a  local version of Zak's theorem.

\proclaim{Question 7} 
 Let $X^n\subset\pp{ \binom{n+2}2-1}$ 
  be a patch of a
complex manifold not contained in a hyperplane. Let $x\in X$ be
a general point and let $v\in T_xX$ be a generic tangent vector. 
If $III^v=0$, must $X$ be a patch of the Veronese $v_2(\pp n)$?
\endproclaim

  The difference between knowing that 
the cubic form is zero and knowing 
that the refined third fundamental
form with respect to all  tangent directions is zero is
a difference of $\binom{n+1}2\binom{n+2}3-(n+1)\binom{n+1}2 $ vs 
$\binom n2\binom n3 -\binom{n+1}2 $ equations. 
From the proof of theorem 3, one sees
how to construct the germ of a negative answer, but there is no
reason to believe any germ will extend to a smooth variety.

\smallpagebreak
 
 Using    proposition 7, we
obtain  a stronger
infinitesimal rigidity result 
than in [L2] by   observing that
if $\ii$ is the second fundamental form of
a Severi variety, then $\ii\up 1=0$
and thus $(II/\tann (v))\up 1=0$.  Theorems 1,2,3 in
the case of Severi varieties follow  from proposition
4 and this observation. The proofs 
given here of theorems 1 and 2
are   better than those in [L2] because here the basis
vectors used for
$T_xX$ are in the baselocus of $\ii$. In contrast, in [L2] a basis
consisting of essentially generic vectors (although not a generic basis)
was used.  One could write out a corresponding better proof for
the $\Bbb O\pp 2$ case as well.

\proclaim{Question 8} Let  $X^n\subset\pp\na$
be a smooth variety with degenerate tangential variety.
  Let $x\in X$ be a general 
point and let $v\in T_xX$ be a generic tangent vector. With the
notation of the paragraphs above consider
$|II_{X,x}|/\tann (v)$ as a system of quadrics on
$\tsingloc (\tann (v))$. Must 
$
 (|II_{X,x}|/\tann (v))\up 1=0?
$
\endproclaim

An affirmative answer to question 8 would
provide a new proof
of the Fulton-Hansen theorem relating the dimensions
of $\sigma (X)$ and $\tau (X)$ that is
differential-geometric and elementary in nature
in the  
case   $X$ is smooth.
(In particular, one would not need 
Deligne's Bertini theorem.)

\smallpagebreak
 
A  variant  of question 8 is as follows:
Let $A\subset S^2\bcc n$ be a system of
quadrics with a tangential defect (i.e. the quadrics
in $A$ satisfy a polynomial relation). What additional
conditions can one impose on $A$ to imply that
if $A=|II_{X,x}|$ where $x\in X$ is a general point, then any
tensor corresponding to $|III^v_{X,x}|$ must be zero? 

\subheading{Ideas towards more geometric proofs}
While the proofs here are rather short, it would 
be desirable to have more geometric arguments. The rigidity statements
in [HM] are proven by exploiting that $\La 2 T_xX$ is generated
by elements of the form $v\ww v'$ where $[v]\in \tbaseloc\ii$
and $v'\in T_{[v]}\tbaseloc\ii$. In [LM] we  show  that
if $X$ is homogeneous, if $\sigma (\tbaseloc\ii )=\ppp T_xX$,
then $\ii\up 1=0$, so in particular $III^v\equiv 0$. Thus,
part of the results here follow from geometric arguments,
but it is not in general true that all third order information
can be recovered from $III^v$. Accordingly, some further geometric
properties of $\tbaseloc\ii$ are needed.

\subheading{Other  open problems} 
Lebrun's program to classify the
quaternionic-K\"ahler
manifolds with postive scalar curvature (see [Le], [LS])
has reduced (via the twistor transform)
the classification problem to classifying the
contact Fano manifolds. (It is generally
conjectured that the only  quaternionic-K\"ahler
manifolds with postive scalar curvature are homogenous.)
The only known contact Fano manifolds are the adjoint
varieties. Given a contact Fano manifold, the base locus
of its second fundamental form must be a Legendrian
variety. So it is 
of particular   interest
 to determine the extent a variety
must resemble a homogeneous Legendrian variety (resp. adjoint variety)
before one can conclude that it is 
 a homogeneous Legendrian variety
(resp. adjoint variety).

Another problem is to determine rigidity for the cases of
$\pp 1\times\pp n$ and $G(2,5)$, which are not covered by
the theorems above. By Ein's results
on dual varieties, one  would conjecture
that these varieties are rigid to second order as well.
If these cases are rigid to order two, it would
give strong evidence for an affirmative answer
to the following question:

\proclaim{Question 9}
Let $G/P\subset\pp N$ be a homogeneous variety 
in its canonical embedding with
$\ff k\neq 0$, $\ff{k+1}=0$. Assume  
  $G/P$ is not   a 
quadric hypersurface. Let
$X^n\subset\pp\na$ be a patch of a complex manifold and let $x\in X$ be a
general point. If $\fff d{X,x}=\fff d{G/P}$ for all $d\leq k$,
must $X$ be an open subset of $G/P\subset\pp N$?
\endproclaim

A weaker version of this question would be to require
additionally that all differential invariants of $X$
other than the fundamental forms are zero up to order $k$. 

\heading Proofs\endheading

We will use formulas for projective
differential invariants derived in [L1].

The idea of the proofs is as follows: given
any variety $X\subset\ppp V$, one has
the first order adapted frame bundle
  $\pi :\cf^1_X\ra X$.
The elements $f\in \cf^1=\cf^1_X$ are bases of $V$ that respect the
flag $\hat x\subset\hat T_xX\subset V$,
where $\hat x$ is the line in $V$ corresponding to $x$
and $\hat T_xX$ is the cone over the embedded
tangent space. In particular,
each $f\in\cf^1$ determines a splitting of the
flag which we denote $\hat x + T + N$.
 Although it is not in general a Lie group,
$\cf^1\subset GL(V)$.

Write the pullback of the Maurer-Cartan form of $Gl(V)$ to
$\cf^1$
as
$$
\Omega = \pmatrix \oo 0&\ooo 0\beta&\ooo 0\nu\\
\oo\alpha&\ooo\alpha\beta&\ooo\alpha\nu\\
0&\ooo\mu\beta&\ooo\mu\nu\endpmatrix
$$
with index ranges $1\leq\alpha,\beta\leq\tdim X$,
$\tdim X+1\leq\mu,\nu\leq \tdim\ppp V$.

If $X=G/P$, one can reduce
$\cf^1$ until it is isomorphic to $G$ (with fiber
isomorphic to $P$). In that case one obtains the
Maurer-Cartan form symbolically as:
$$
\Omega_G = \pmatrix \oo 0&\ooo 0\beta&0\\
\oo\alpha&\ooo\alpha\beta=\rho_T(\frak h) &\ooo\alpha\nu
=A_2(\ooo 0\beta ) \\
0&\ooo\mu\beta=
A_1(\oo \alpha ) &\ooo\mu\nu=\rho_N(\frak h) \endpmatrix
$$
where 
$H$ is the semi-simple part of $P$,
$T=T_xX$, $N=N_xX$ are $H$-modules,
and
$A_1,A_2$ are $H$-equivariant maps.
  The zero in the upper
right hand block indicates that  any infinitesimal change in
the splitting statisfies the
\lq\lq transversality\rq\rq \
condition that $dN\subseteq \{ T + N\}$.
The dependence of the
$\ooo\alpha\nu$ block   on the forms $\ooo 0\beta$
  indicates that if one changes the choice of $T$,
there is a corresponding change in choice of $N$ mandated.

If $X$ is a variety with the same second fundamental form as $G/P$,
by restricting bases in $T_xX$ and $N_xX$ to be
we can reduce $\cf^1_X$ to a bundle $\cf^2_X$
where the pullback of the the Maurer-Cartan
form looks like:
$$
\Omega  = \pmatrix \oo 0&\ooo 0\beta&\ooo 0\nu\\
\oo\alpha&\ooo\alpha\beta=\rho_T(\frak h)+ w_1 &\ooo\alpha\nu
  \\
0&\ooo\mu\beta  &\ooo\mu\nu=\rho_N(\frak h)+ w_2 \endpmatrix
$$
where $w_1,w_2$ are linear combinations of the other
forms appearing in the Maurer-Cartan form.
The proofs proceed  by showing that there are  
reductions of $\cf^2_X$ to $G$  by restricting the admissible
splittings that reduce to
$\Omega_G$.

In practice, we prove  this by showing the
  invariants $F_k\in \pi^*(S^kT^*X\ot NX)$  that contain
the geometric information of the $k$-th derivative of the 
embedding $X\ra \pp N$, are zero for $k>2$.
 In
frames one writes
$F_k=\rr\mu{\alpha_1}\hdots{\alpha_k}\oo {\alpha_1}\hdots\oo
{\alpha_k}\ot\ee\mu$, where $\oo\beta$ is a basis of the semi-basic forms
and $\ee\mu$ is a basis of $N_xX(1)$, and the
$\rr\mu{\alpha_1}\hdots{\alpha_k}$ are functions defined on $\cf^1$.
$F_k$ measures the infinitesmal motion of $X$ away from its embedded
tangent space to $(k-1)$-st order.

We recall the following formulae from [L1]:
$$
\align
\rr\mu\alpha\beta\gamma\oo\gamma &=
-d\qq\mu\alpha\beta - \qq\mu\alpha\beta\ooo 00 -\qq\nu\alpha\beta\ooo\mu\nu
+\qq\mu\alpha\delta\ooo\delta\beta + \qq\mu\beta\delta\ooo\delta\alpha
\tag L1 2.15\\
 \rr\mu\alpha\beta{\gamma\delta}\oo\delta &=
-d\rr\mu\alpha\beta\gamma -2\rr\mu\alpha\beta\gamma\ooo 0 0
-\rr\nu\alpha\beta\gamma\ooo\mu\nu + \frak
S_{\alpha\beta\gamma}\rr\mu\alpha\beta\delta\ooo\delta\gamma\tag 
L1 2.17 \\
& \ \ \
-\frak S_{\alpha\beta\gamma}
\qq\mu\alpha\delta\qq\nu\beta\gamma\ooo\delta\nu +
\frak S_{\alpha\beta\gamma}\qq\mu\alpha\beta\ooo 0\gamma  \\
 \rr\mu\alpha\beta{\gamma\delta\epsilon}\oo\epsilon &=
-d\rr\mu\alpha\beta{\gamma\delta} 
-3\rr\mu\alpha\beta{\gamma\delta}\ooo 0 0
-\rr\nu\alpha\beta{\gamma\delta}\ooo\mu\nu \tag 
L1 2.20\\
& \ \ \ +  \frak
S_{\alpha\beta{\gamma\delta}}(
\rr\mu\alpha\beta{\gamma\epsilon}\ooo\epsilon\delta
+2\rr\mu\alpha\beta{\gamma }\ooo 0\delta\\
& \ \ \ 
-(\rr\mu\alpha\beta\epsilon\qq\nu\gamma\delta
+\rr\nu\alpha\beta\gamma\qq\mu\delta\epsilon)\ooo\epsilon\nu
-\qq\mu\alpha\beta\qq\nu\gamma\delta\ooo 0\nu )
\endalign
$$
The functions $\rr\mu{\alpha_1}\hdots{\alpha_k}$   vary in the fiber as
follows: Under a motion
$$
\align
\ee\alpha&\mapsto \ee\alpha + g^0_{\alpha}\ee 0\\
\ee\mu&\mapsto\ee\mu +g^{\alpha}_{\mu}\ee\alpha + g^0_{\mu}\ee 0
\endalign
$$
the corresponding changes in the coefficients of $F_3,F_4$ are
as follows:
$$
\align
&\Delta\rr\mu\alpha\beta\gamma =  
\frak S_{\alpha\beta\gamma}g^0_{\alpha}\qq\mu\beta\gamma +
\frak
S_{\alpha\beta\gamma}g^{\delta}_{\nu}\qq\nu\alpha\beta\qq\mu\gamma\delta
\tag L1 2.24\\
&\Delta\rr\mu\alpha\beta{\gamma\delta} =
\frak S_{\alpha\beta\gamma\delta}g^0_{\alpha}\rr\mu\beta\gamma\delta +
\frak S_{\alpha\beta\gamma\delta}g^{\epsilon}_{\nu}
(\rr\nu\alpha\beta\gamma\qq\mu\delta\epsilon+
\qq\nu\alpha\beta\rr\mu\gamma\delta\epsilon)
+ g^0_{\nu}\qq\mu\alpha\beta\qq\nu\gamma\delta .
\endalign
$$
\demo{Proof of theorem 1}
$z=0$ because the prolongation of $\ii$  is zero.
Let $V$ have basis $\{ \ee 0,\ee i,\ee s,\dee si\}$,
$1\leq i,j,k,l\leq n$, $n+1\leq s,t,u,v\leq n+m$
adapted such that $x=[\ee 0]$,
$\hat T_xX=\{\ee 0,\ee i,\ee s\}$  and
$
II_{X,x}= \oo i\oo s\ot\dee is
$
(see e.g. [GH], or [L1]).
Note that   the forms $\ooo ij,\ooo st$
are all independent and independent of the semi-basic forms because
they represent infinitesimal motions that preserve our normalization
of $II$.

To show $X$ is a patch of the Segre, we will show all higher
differential invariants are zero.
We see immediately that
$$
\align
&
\rr{si}tu\beta =0 \ \forall \beta\text{ and }t, u\neq s \\
&\rr{si}jk\beta = 0 \ \forall \beta\text{ and }j, k\neq i \endalign
$$
(these equations include the equations for the refined third fundamental
form being zero).
The   nonzero coefficients of $F_3$ satisfy the
following equations.
(Here and in what follows, we use the convention that if an index
appears more than twice it is {\it not} to be summed over.
E.g. there is no sum over $s$ in (s1).)
From now on, if latin indicies are distinct, we assume they
are not equal.
$$
\align
&\rr{si}sts\oo s + \rr{si}stk\oo k+ \rr{si}sti\oo i = \ooo it
\tag s1\\
&\rr{si}ijs\oo s +\rr{si}ijt\oo t + \rr{si}iji\oo i = \ooo sj\tag s2\\
&\rr{si}ss\beta\oo \beta = 2\ooo is\tag s3\\
&\rr{si}ii\beta\oo \beta = 2\ooo si\tag s4
\endalign
$$
Since the right hand side of 
(s1), resp. (s2), is independent of 
$s$, resp. $i$, we conclude (assuming $n,m\geq 2$)
$ \rr{si}sts=0 ,
 \rr{si}iji=0$ and
$$
\align
&\rr {si}stk=\rr{ui}utk \tag s5\\
&\rr {si}ijt=\rr{sk}kjt 
\tag s6\endalign
$$
Now
$$
\align
&\Delta \rr {si}stj=g^i_{(tj)}\\
&\Delta \rr {si}ijt=g^s_{(jt)}\\
&\Delta \rr {si}sti=g^i_{(ti)}+g^0_t\\
&\Delta \rr {si}ijs=g^s_{(js)}+g^0_j\endalign
$$
Fixing a
particular $(i,s)$,
use $g^i_{tj},g^s_{tj}, g^i_{ti},
g^s_{js}$ to normalize all these terms to zero. By (s5,s6) the
normalizations send the terms to zero for all $i,s$.
Now (s1,s2) imply $\ooo it=0,\ooo ti=0$ for all $i,t$ 
so (s3,s4) imply
$$
\align
&\rr {si}ss\beta =0\\
&\rr {si}ii\beta =0.\endalign
$$

We have now reduced to frames where $F_3=0$. (At this point one has a
projective connection on $TX$ isomorphic to that on the Segre.)
The coefficients of
$F_4$, $\rr{si}\alpha\beta{\gamma\delta} $ are zero
unless among $\alpha\beta\gamma\delta $,
two are in the $s,t,u$ range and two are in the $i,j,k$ range,
and at least one is equal to $s$ and one equal to $i$.
We   use  
$$
\Delta\rr {si}ss{ii}= g^0_{si}
$$
to normalize $\rr {si}ss{ii}  =0$ for all $s,i$.
 This uses up all
the freedom to normalize differential invariants. (The
$g^0_{\alpha}$ terms were useless as they always appeared with
a corresponding $g^{\epsilon}_{\mu}$ term.)

The remaining coefficients of $F_4$ that are potentially
nonzero are $\rr {si}st{ji},\rr{si}ss{ij},\rr {si}ii{st}$.


We now  examine the coefficients of $F_5$. Fortunately
most of these are immediately seen to be zero.  
$\rr{si}\alpha\beta{\gamma\delta\epsilon}$ is zero if
four or all of the lower indices are all in the 
same range. Moreover there must be at least two indicies
that are either $i$ or $s$. Consider
$$
\align
\rr{si}su{jks}\oo s+\rr{si}su{jki}\oo i
&=\rr{si}su{ji}\ooo ik + \rr{si}su{ki}\ooo ij\\
\rr{si}tu{iks}\oo s+\rr{si}tu{iki}\oo i
&=\rr{si}su{ik}\ooo st + \rr{si}st{ik}\ooo su\\
 \rr{si}ss{jki}\oo i
&=\rr{si}su{ji}\ooo ik + \rr{si}su{ki}\ooo ij\\
\rr{si}tu{iks}\oo s+\rr{si}tu{iki}\oo i
&=\rr{si}su{ik}\ooo st + \rr{si}st{ik}\ooo su .\endalign
$$
In all equations the forms on the right hand side are all independent
and independent of the forms on the left hand side (which are
independent as well). Thus all coefficients appearing are zero,
in particular, all coefficients of $F_4$ are zero.
Now consider
$$
\rr{si}st{jii}\oo i +\rr{si}st{jis}\oo s=\ooo 0{tj}
$$
Since the right hand side is independent of $i,s$, we conclude
both sides are zero. Using that $\ooo 0{tj}=0$
for all $t,j$,
  the equations   $\rr{si}ss{ii\beta}\oo\beta= 2\ooo
0{si}$ 
 imply $F_5=0$.
 We easily see the coefficients of $F_6$ are all zero and
thus all higher differential invariants are zero.\qed\enddemo

\demo{Proof of theorem 2}
Again $z=0$ because $\ii\up 1=0$.
 
Let $V$ have basis $\{\ee 0,\dee 1j,\dee 2j,\dee jk\}$, where 
$3\leq j,k,l\leq n+2$, $\{\alpha\}=\{ 1j, 2j\}$.
Normalize such that
$
II= (\oo{1j}\oo{2k}-\oo{1k}\oo{2j})\ot\dee jk, \ j<k
$.
Note that   the forms $\ooo {1i}{1j},\ooo {2i}{2j},
\ooo {1i}{1i},\ooo {2i}{2i},\ooo {2i}{1i},
\ooo {1i}{2i} $
are all independent and independent of the semi-basic forms because
they represent infinitesimal motions that preserve our normalization
of $II$.
We have
$$
\align &
\rr{ij}{(1k)}{(1l)}\beta = 0\ \forall i,j,k,l\text{ distinct  and }
\forall\beta\tag g1\\
 &
\rr{ij}{(2k)}{(2l)}\beta = 0\ \forall i,j,k,l\text{ distinct  and }
\forall\beta\tag g2\endalign
$$
(these equations include that the refined third fundamental
form is zero).
From now on, assume all indices are distinct.
Using (g1), (g2), we have
$$
\align &
\rr{(ij)}{(1i)}{(1k)}{(1i)}\oo{1i}+
\rr{(ij)}{(1i)}{(1k)}{(1j)}\oo{1j}+
\rr{(ij)}{(1i)}{(1k)}{(2i)}\oo{2i}+
\rr{(ij)}{(1i)}{(1k)}{(2j)}\oo{2j}+
\rr{(ij)}{(1i)}{(1k)}{(2l)}\oo{2l}\tag g3\\
&\ \ \ =\ooo{(2j)}{(1k)} \endalign
$$
 The right hand side of (g3) is independent of $i$,
so comparing with the same expression using $m$ instead of $i$,
(here we use that $n\geq 4$)
we obtain:
$$
\align
&\rr{(ij)}{(1i)}{(1k)}{(1i)}=0\tag g4\\
&\rr{(ij)}{(1i)}{(1k)}{(1j)}=\rr{(mj)}{(1m)}{(1k)}{(1j)}
  \tag g5\\
&\rr{(ij)}{(1i)}{(1k)}{(2i)}=\rr{(mj)}{(1m)}{(1k)}{(2i)}
  \tag g6\\
&\rr{(ij)}{(1i)}{(1k)}{(2j)}=\rr{(mj)}{(1m)}{(1k)}{(2j)}
  \tag g7\\
&\rr{(ij)}{(1i)}{(1k)}{(2l)}=\rr{(mj)}{(1m)}{(1k)}{(2l)}
  \tag g8\endalign
$$
Now
$$
\align
&
\Delta \rr{(ij)}{(1i)}{(1k)}{(2l)}= g^{(2j)}_{(kl)}
\\
&
\Delta \rr{(ij)}{(1i)}{(1k)}{(2j)}= g^{(2j)}_{(jk)}
+g^{0}_{(1k)}
\endalign
$$
Using these equations and the corresponding
equations with the role of $1$ and $2$ reversed, we reduce to frames where
$ 
 \rr{(ij)}{(1i)}{(1k)}{(2l)}= 0, 
  \rr{(ij)}{(1i)}{(1k)}{(2j)}=0,
 \rr{(ij)}{(2i)}{(2k)}{(1l)}= 0,
  \rr{(ij)}{(2i)}{(2k)}{(1j)}=0$. 
In these frames, $\ooo{1k}{2j},\ooo{2k}{1j}=0$ hence
$$
0=\rr {(ij)}{(1i)}{(1i)}{\beta}\oo{\beta}=-2\ooo{2j}{1i}\tag g9
$$
and similarly with the role of $1$ and $2$ reversed.
Thus the only nonzero terms left in $F_3$ are
$ \rr{(ij)}{(1i)}{(1j)}{(1k)},\rr{(ij)}{(2i)}{(2j)}{(2k)}$.
Consider
$$
\rr {(ij)}{(1i)}{(1j)}{(1i)\beta}\oo{\beta}=2
 \rr{(ij)}{(1i)}{(1j)}{(1k)}\ooo{1k}{1i}\tag g10
$$
Both sides of (g10) must be zero because the forms
$\ooo{1k}{1i}$ are all independent and independent of
the semi-basic forms. The analogous equation holds
with $2$'s. Hence we see $F_3=0$.

To have a nonzero coefficient of $F_4$, $\rr{(ij)}\alpha
\beta{\gamma\delta}$, in the lower indicies there must
be two $1$'s and two $2$'s, and at least two of the $k$-indices
must be $i$ or $j$. Consider
$$
\align
&\rr {(ij)}{(1i)}{(1k)}{(2l)(2j)}\oo{2j}=\ooo{2j}{kl}\tag g11\\
&\rr {(ij)}{(1i)}{(1i)}{(2l)(2j)}\oo{2j}=2\ooo{2j}{il}.\tag g12\endalign
$$
Since the right hand side of (g11) is independent of
$i$, we conclude (after switching the roles of $i$ and $j$)
that $\rr {(ij)}{(1i)}{(1k)}{(2l)(2j)}$ is independent of $i,j$
(with neither $k,l$ equal to $i$ or $j$, but $k=l$ is possible). Using
$$
\Delta \rr {(ij)}{(1i)}{(1k)}{(2l)(2j)}= g^0_{kl}
$$
we normalize all these terms to zero. This implies $\ooo{2j}{il}=0
$
and thus $\rr {(ij)}{(1i)}{(1i)}{(2l)(2j)}=0$ for all $i,j,l$ distinct
as well, and similarly with the role of $1$ and $2$ reversed.
Thus the remaining nonzero terms in $F_4$ are
$\rr {(ij)}{(1i)}{(1i)}{(2j)(2j)},\rr {(ij)}{(1i)}{(2i)}{(1j)(2j)}$.
Consider
$$
\align 
&\rr {(ij)}{(1i)}{(2j)}{(1k)(2l)(1i)}\oo{1i}+
\rr {(ij)}{(1i)}{(2j)}{(1k)(2l)(1j)}\oo{1j}\tag g13\\
&+
\rr {(ij)}{(1i)}{(2j)}{(1k)(2l)(2i)}\oo{2i}+
\rr {(ij)}{(1i)}{(2j)}{(1k)(2l)(2j)}\oo{2j}\\ 
& = -\ooo 0{kl}.\endalign
$$
Since the right hand side of (g13) is independent of $i,j$, we
conclude $\ooo 0{kl}=0$ and hence the left hand side is
zero as well. Now it is easy to see the rest of the terms in
$F_5$ are zero and all higher forms are zero.\qed
\enddemo

\subheading{Remark} While the last step used $m\geq 4$
a second time, there is an
alternate argument that avoids it here.
One first observes the  $\ooo 0{kl}$
are semi-basic and then uses the equation for
$\rr {(ij)}{(1i)}{(1i)}{(2j)(2k)\beta}\oo{\beta}$. 
 
\demo{Proof of theorem 4}
 $III=0$ and $X$ not contained in a hyperplane implies
that $\ii=\ppp S^2T^*$. Take a basis
$(\ee 0,\ee j\ee{ij})$ of $V$, $1\leq i,j,k\leq n$ such that
$II= \oo i\oo j\ot\ee{ij}$. With this normalization
the forms
$\ooo ij$ are all independent and independent of the
semi-basic forms.
We cannot use  the $g^i_{jk},g^0_j$
to make normalizations
because we
  assumed $F_3=0$. The coefficients of 
$F_4$, $\rr{ij}kl{mp}$, must be zero if three or four of the
lower coeffients are different from $i,j$. Regarding
the other coefficients,
$$
\align
\rr{ij}ik{li}\oo i + \rr{ij}ik{lj}\oo j
&= \ooo j{kl}\tag v1\\
\rr{jj}jk{lj}\oo j  
&= 2\ooo j{kl}\tag v2\endalign
$$
and one has the corresponding equations with $k=l$.
Combining (v1) and (v2) we conclude $\rr{ij}ik{li},\rr{ij}ik{ki}=0$,
(here we use $n\geq 3$, see [GH] or [L2] for a proof when $n=2$) and
if $n\geq 4$ we also have,
since the right hand side of (v1) is independent of $i$,
$
\rr{ij}ik{lj}=\rr{mj}mk{lj} 
$.
The   variability
$$
\align
&\Delta\rr{ij}ik{lj} =g^0_{kl}\\
&\Delta\rr{ij}ik{kj} =g^0_{kk}\endalign
$$
  allows us to normalize
$ \rr{ij}ik{lj},\rr{ij}ik{kj}=0$ which implies $\ooo i{kl},
\ooo i{kk}=0$, in turn
implying $\rr{ii}ik{li},\rr{ii}ik{li} =0$.
Since $\ooo j{ik},\ooo j{ii}=0$ we have
$$
\align 
&\rr{ij}ii{k\beta}\oo \beta  
  =0\\
&\rr{ii}ii{k\beta}\oo \beta  
  =0\endalign
$$
so $\rr{ij}ik{ii},\rr{ij}ii{kj},\rr{ij}ii{ii},\rr{ij}ii{ij}=0$.
The only potentially nonzero terms in $F_4$ are
$\rr{ij}ij{ij},\rr{ii}ii{ii}$. Consider
$$
\rr{ij}ij{kl\beta}\oo \beta = \ooo 0{kl}
$$
which implies $\ooo 0{kl}$ is semi-basic. Using
$$
\rr{ij}ii{jk\beta}\oo \beta = \rr{ij}ii{jj}\ooo jk +
\ooo 0{ik}, 
$$
we conclude $\rr{ij}ii{jj}=0$,
since $\ooo jk$ is independent of the semi-basic forms.
Using the corresponding
equation with $i$ replacing $j$ we see $F_4=0$. Now
it is easy to see that $F_5$ and all higher invariants
are zero. \qed\enddemo

\subheading{Acknowledgements} It is a pleasure to thank
Jun-Muk Hwang for useful conversations.

\enddocument

this is unnecessary:

 For example
$$
\align 
\rr{ij}ij{klj}\oo j &=  
-\ooo 0{kl}\tag v4\\
\rr{ii}ii{kli}\oo i &=  
+\ooo 0{kl}\tag v5\endalign
$$
 And one easily sees all coefficients are zero since the right
hand side is independent of $j$.
\enddemo

\enddocument